\documentclass[a4paper,11pt]{article}
\usepackage{paralist}   

\usepackage{amsmath}
\usepackage{amssymb}
\usepackage{amsfonts}
\usepackage{color}

\usepackage{graphicx}
\usepackage{subfig}
\usepackage{bm}
\usepackage{tikz}

\newtheorem{theorem}{Theorem}
\newtheorem{definition}{Definition}
\newtheorem{proposition}[theorem]{Proposition}
\newtheorem{corollary}[theorem]{Corollary}
\newtheorem{lemma}[theorem]{Lemma}
\newtheorem{remark}[theorem]{Remark}
\newtheorem{example}{Example}

\begin{document}

\title{\textbf{The continuous nonstationary Gabor transform on LCA groups with applications to representations of
the affine Weyl-Heisenberg group}}
\author{Michael Speckbacher and Peter Balazs\footnote{Acoustics Research Institute, Austrian Academy of Sciences,  Wohllebengasse 12-14, 
1040 Vienna, Austria, speckbacher@kfs.oeaw.ac.at, peter.balazs@oeaw.ac.at}}
\date{}
\maketitle

\begin{abstract}
In this paper we introduce and investigate the concept of reproducing pairs which generalizes  continuous frames.
We will introduce a concept that represents a unifying way to look at certain continuous frames (resp. reproducing pairs) on LCA groups,
which can be described as 
continuous nonstationary Gabor systems and investigate conditions for these systems to form a  continuous frame (resp. reproducing pair).
As a byproduct we identify the structure of the frame operator (resp. resolution operator).
Moreover, we ask the question, whether there always exist mutually dual systems with the same structure such that the resolution operator
is given by the identity, i.e. given $A:X\rightarrow B(\mathcal{H})$, if there exist $\psi,\varphi\in\mathcal{H}$, s.t.
\begin{equation*}
 f=\int_X \langle f,A(x)\psi\rangle A(x)\varphi d\mu(x),\ \ \forall f\in \mathcal{H}
\end{equation*}
and show that the answer is not affirmative. As a counterexample we use a system generated by a unitary action of a subset of 
the affine Weyl-Heisenberg group in $L^2(\mathbb{R})$.
\end{abstract}
\vspace{1cm}
\textbf{Math Subject Classification:} 22B99, 43A32, 42C15, 42C40.
\textbf{Submitted to:} Journal of Physics A: Mathematical and Theoretical

\newpage

\section{Introduction}
Motivated by physical applications \cite{alanga93a,ka90}, in order to generalize the coherent states approach, the concept of continuous 
frames has been introduced in the early 1990's independently by  Ali et al. \cite{alanga93} and Kaiser \cite{ka94}.
Coherent states, see e.g. \cite{alanga00}, are widely used in many areas of theoretical physics, in particular in quantum mechanics, where 
classical coherent states are generated by a group action on a single mother wavelet and lead to a resolution of the identity. More general,
continuous frames yield a resolution of a positive, bounded and invertible operator. 

This raises the question if such a resolution necessitates the frame property. To answer this issue, we will introduce the idea of a 
reproducing pair where a pair of mappings, in place of a single mapping, is used for an invertible analysis/ synthesis process and compare 
reproducing pairs to continuous frames.
For discrete frames this question of dual systems is a current topic of research, although this is bounded in most cases 
to the frame property \cite{raey}, and often to certain types of frames, see e.g. \cite{persampta13,Christensen2014198}. Introducing 
reproducing pairs, we ask the more general question of whether reconstruction is possible, without assuming the frame property a-priori. 
This is related to the topic of frame multipliers \cite{xxlmult1,xxlbayasg11}. These are operators consisting of analysis, element-wise
multiplication 
with a fixed symbol, and synthesis, and appear in a lot of scientific disciplines. In Physics they represent the link between classical 
and quantum mechanics, so called quantization operators \cite{alanga00}.
The invertibility of multipliers is a central topic \cite{stobalrep11,balsto09new,uncconv2011} in the mathematical investigation of
these operators. This includes the question, when a system of two mappings forms a reproducing pair.

Tight frames, i.e. systems where the corresponding frame operator is a multiple of the identity, are often preferred to non-tight 
ones since the inversion
of the frame operator is straightforward. This begs the question
if, given a particular structure, there always exist a tight frame, or reproducing pair with resolution operator $\lambda I$, with this 
structure. This will be answered in this paper for the particular class of continuous systems investigated here.

Two of the most widely used continuous frame transforms are the Short Time Fourier Transform (STFT) \cite{groe1} and the Wavelet transform 
\cite{ma09}, in particular in signal processing and acoustics. In those applications they are used in their sampled, discretized version. 
Both transforms have a time-frequency resolution, which is either fixed for all frequencies (for the STFT), or follow a given rule 
(for wavelets). 
In practice functions or signals often show particular time-frequency characteristics which call for adaptive and adaptable representation 
\cite{badokowto13}. In \cite{nsdgt10} the authors introduced adaptivity either in time or frequency, where perfect reconstruction is still
possible. This ansatz will be adapted within this paper to introduce continuous nonstationary Gabor frames. We will 
see in Section \ref{cont-nonstationary gabor} that the transforms considered here can be rewritten in terms of a convolution and recall 
that the Fourier transform diagonalizes convolutions (within the right function spaces). Hence it seems natural to consider locally 
compact abelian (LCA) groups as domain of definition, as Fourier theory provides sufficient results on these groups and the Hilbert space 
$\mathcal{H}=L^2(G)$. Within this setting we will then derive a sufficient substitute for the frame (resp. reproducing pair) condition 
independent of $f\in L^2(G)$ and show that the frame (resp. resolution) operator is given by a Fourier multiplier.

The study of representations of the affine Weyl-Heisenberg group is of particular interest since it 
contains both the Weyl-Heisenberg group and the affine group, the underlying groups of the short-time Fourier transform the continuous 
wavelet transform. Recently a lot of effort has been put in the study of intermediate transforms as the $\alpha$-transform, see 
\cite{daforastte08,dalo06}. The continuous nonstationary Gabor transform is a more general concept, also including this setting,
giving rise to a wide range of time frequency transforms. This will give a counterexample to the question
whether there always exist a reproducing pair with a given structure such that the resolution operator is the identity.

The present paper is organized as follows: in Section \ref{prel} we will briefly present the basic results on Fourier analysis on LCA
groups and continuous frames and introduce 
the concept of reproducing pairs. Section \ref{cont-nonstationary gabor} is concerned with the continuous nonstationary Gabor 
transform on LCA groups. The results are then applied in Section \ref{aWH-section}
to representations of subsets of the affine Weyl-Heisenberg group and we will show particularities of reproducing 
pairs in comparison to continuous frames with the help of a particular example.

\section{Preliminaries}\label{prel}
For a short and self-contained introduction to continuous frames, see \cite{ranade06}.
\begin{definition}\label{def-cont-frame}
Let $\mathcal{H}$ be a Hilbert space and $(X,\mu)$ be a measure space.
A mapping $\Psi:X\rightarrow \mathcal{H}$  is called a continuous frame if
\begin{enumerate}[(i)]
 \item $\Psi$ is weakly measurable, i.e. $x\mapsto\langle f,\Psi(x)\rangle$ is a measurable function for all $f\in\mathcal{H}$ 
 \item there exist positive constants $A,B>0$ s.t. 
 \begin{equation}\label{frame-condition}
  A\left\|f\right\|_\mathcal{H}^2\leq\int_{X}\left|\langle f,\Psi(x)\rangle\right|^2d\mu(x)\leq B\left\|f\right\|_\mathcal{H}^2,\ \
  \forall f\in\mathcal{H}
 \end{equation}
\end{enumerate}
\end{definition}
The mapping $\Psi$ is called Bessel if the second, upper inequality in (\ref{frame-condition}) is satisfied.

The standard setting for frames is the discrete version, see for example \cite{christ1}, which is a specialization of this definition
and can be reached by choosing $X$ to be a countable set and $\mu$ the counting measure. Throughout this paper we will assume that 
$\Psi$ is uniformly bounded, i.e. $\sup_{x\in X}\|\Psi(x)\|_\mathcal{H}\leq C$.

Let us define the basic operators in frame theory, the analysis operator
\begin{equation*}
 V_\Psi:\mathcal{H}\rightarrow L^2(X,\mu),\ \ \ V_\Psi f(x):=\langle f,\Psi(x)\rangle
\end{equation*}
and its adjoint operator,  called synthesis operator
\begin{equation*}
 V^\ast_\Psi:L^2(X,\mu)\rightarrow \mathcal{H},\ \ \ V^\ast_\Psi \varphi:=\int_X\varphi(x)\Psi(x)d\mu(x)
\end{equation*}
where the integral is to be understood in the weak sense.
By composition of $V_\Psi$ and $V^\ast_\Psi$ we obtain the frame operator
\begin{equation*}
 S_\Psi:\mathcal{H}\rightarrow \mathcal{H},\ \ \ S_\Psi f:=V^\ast_\Psi V_\Psi f=\int_X\langle f,\Psi(x)\rangle \Psi(x)d\mu(x)
\end{equation*}
The frame operator is obviously self-adjoint and the frame bounds guarantee that it is positive, bounded and invertible. 
The mapping $S^{-1}_\Psi\Psi$ is also a frame, called the canonical dual frame, with frame bounds $B^{-1},A^{-1}$.

The analysis operator $V_\Psi$ is in general not onto $L^2(X,\mu)$ but satisfies a reproducing kernel equation: for $F\in L^2(X,\mu)$, 
there exists $f\in\mathcal{H}$, s.t. $F(x)=V_\Psi f(x)$ if and only if
$F(x)=R(F)(x)$ where $R$ is an integral operator with kernel $\mathcal{R}(x,y):=\langle S_\Psi^{-1} \Psi(y),\Psi(x)\rangle$ and
\begin{equation*}
 R(F)(x):=\int_X\mathcal{R}(x,y)F(y)d\mu(y)
\end{equation*}
 $R$ is moreover the orthogonal projection operator onto the image of $V_\Psi$.

\begin{theorem}
 Let $\Psi$ be a frame, then the following inversion formula holds
 \begin{equation}\label{inversion-general-frame}
  f=\int_{X}\langle f,\Psi(x)\rangle S^{-1}_\Psi \Psi(x)d\mu(x)=\int_{X}\langle f, S^{-1}_\Psi \Psi(x)\rangle \Psi(x)d\mu(x),\ \forall
  f\in\mathcal{H}
 \end{equation}
\end{theorem}
If $\Psi_d$ is another frame satisfying
\begin{equation*}
 f=\int_{X}\langle f,\Psi(x)\rangle \Psi_d(x)d\mu(x)=\int_{X}\langle f,\Psi_d(x)\rangle \Psi(x) d\mu(x),\  \ \forall f\in\mathcal{H}
\end{equation*}
then $\Psi_d$ is called a dual frame. In most cases there exist several dual frames since in general ker$(V^\ast_\Psi)\neq \{0\}$.

If we do not restrict to an expansion via mutually dual frames we can introduce the following definition motivated by 
\cite[Definition 1.1.1.33]{te01}. To this end, we denote the space of bounded linear operators
 with bounded inverse from $\mathcal{H}$ to $\mathcal{K}$ by  $GL(\mathcal{H},\mathcal{K})$. 
\begin{definition}\label{rep-pair-definition}
 Let $(X,\mu)$ be a measure space and $\Psi,\Phi:X\rightarrow\mathcal{H}$ weakly measurable.
 The pair of mappings $(\Psi,\Phi)$ is called a reproducing pair for $\mathcal{H}$ if the resolution operator 
 $C_{\Psi,\Phi}:\mathcal{H}\rightarrow \mathcal{H}$, weakly defined by
 \begin{equation}\label{rep-pair-def}
 C_{\Psi,\Phi} f:=\int_X \langle f,\Psi(x)\rangle \Phi(x)d\mu(x)
 \end{equation}
is an element of $GL(\mathcal{H})$.
\end{definition}
Note that this definition is indeed a generalization of continuous frames because if we choose $\Psi=\Phi$, then 
the above definition of reproducing pair corresponds to the original definition of a continuous frame in \cite{alanga93}.

In the same paper the authors showed that it is possible to generate new continuous frames  ``equivalent'' to a given continuous frame. 
In this paper, we will adapt a more general concept from \cite{jale14} to show the result for reproducing pairs. 

 \begin{lemma}\label{reprod-equiv-lemma}
Let $\mathcal{K}$ be a Hilbert space, $(Y,\mu')$ be a measure space, $\rho:Y\rightarrow X$ a bijective mapping which satisfies 
$\mu'\circ \rho^{-1}=\mu$ and preserves measurability,
 $T\in GL(\mathcal{H},\mathcal{K})$ and $\tau:Y\rightarrow \mathbb{C}$ a measurable function with $|\tau(y)|=1$.\\
  If we define $\widetilde \Psi(y):=\tau(y) T ( \Psi\circ \rho) (y)$ (and $\widetilde \Phi$ respectively), 
  then $(\Psi,\Phi)$ is a reproducing pair for $\mathcal{H}$ with respect to $(X,\mu)$, if and only if $(\widetilde \Psi,\widetilde \Phi)$ 
  is a reproducing pair for $\mathcal{K}$ with respect to $(Y,\mu')$.
 \end{lemma}
 \textbf{Proof:} Let $f,g\in\mathcal{K}$. It holds
 \begin{align*}
  \langle C_{\widetilde \Psi,\widetilde \Phi} f,g\rangle_{\mathcal{K}} &=\int_{Y}\big\langle f,\tau(y)T\Psi(\rho(y))\big\rangle_{\mathcal{K}}
 \big \langle \tau(y)T\Phi(\rho(y)),g\big\rangle_{\mathcal{K}} d\mu'(y)\\
&=\int_{Y}\big\langle T^\ast f,\Psi(\rho(y))\big\rangle_{\mathcal{H}}\big\langle \Phi(\rho(y)),T^\ast g\big\rangle_{\mathcal{H}} d\mu'(y)\\
&=\int_{X}\big\langle T^\ast f,\Psi(x)\big\rangle_{\mathcal{H}}\big\langle \Phi(x),T^\ast g\big\rangle_{\mathcal{H}} d\mu(x)\\
&=\langle C_{ \Psi,\Phi} T^\ast f,T^\ast g\rangle_{\mathcal{H}}\\
&=\langle TC_{ \Psi,\Phi} T^\ast f, g\rangle_{\mathcal{K}}
\end{align*}
 Hence, we can identify the resolution operator $C_{\widetilde \Psi,\widetilde \Phi} = TC_{ \Psi,\Phi} T^\ast$ and the result follows as 
 $C_{\widetilde \Psi,\widetilde \Phi}\in GL(\mathcal{K})$ if and only if
 $C_{ \Psi,\Phi}\in GL(\mathcal{H})$.\hfill $\Box$

\vspace{0.3cm}
Unlike the frame operator $S_\Psi$, $C_{\Psi,\Phi}$ is in general neither positive nor self-adjoint, since $C_{\Psi,\Phi}^\ast=
C_{\Phi,\Psi}$. From (\ref{rep-pair-def}) we can also derive a necessary condition for $F:X\rightarrow \mathbb{C}$ to be an element of the 
image of $V_\Psi$ in terms of a reproducing kernel. To do so, we need to define the domain of $V_\Phi^\ast$ as follows 
\begin{equation*}
dom(V_\Phi^\ast):=\Big\{F\in L^\infty(X,\mu):\ \int_X F(x)\Phi(x)d\mu(x)\ converges\ weakly\Big\}
\end{equation*}

\begin{proposition} \label{sec:repordkern1}
 Let $(\Psi,\Phi)$ be a reproducing pair for $\mathcal{H}$ and $F\in dom(V_\Phi^\ast)$. It holds that $F(x)=\langle f,\Psi(x)\rangle$, 
 for some $f\in\mathcal{H}$, if and only if $F(x)=R(F)(x)$ with the integral kernel
 \begin{center}
 $\mathcal{R}(x,y)=\langle C_{\Psi,\Phi}^{-1}\Phi(y),\Psi(x)\rangle$
 \end{center}
 Moreover, $L^1(X,\mu)\cap L^\infty(X,\mu)\subset dom(V_\Phi^\ast)$, which in particular implies that $dom(V_\Phi^\ast)\cap L^2(X,\mu)$ is 
 dense in $L^2(X,\mu)$.
\end{proposition}
\textbf{Proof:}
Let $F(x)=\langle f,\Psi(x)\rangle$, then 
\begin{align*}
 R(F)(x) &=\int_X \langle f,\Psi(y)\rangle\langle \Phi(y),(C_{\Psi,\Phi}^{-1})^\ast\Psi(x)\rangle d\mu(y)\\
&=\langle C_{\Psi,\Phi}f, (C_{\Psi,\Phi}^{-1})^\ast\Psi(x)\rangle=V_\Psi f(x)=F(x)
\end{align*}
Assume now that $R(F)(x)=F(x)$. Since, $F\in dom(V_\Phi^\ast)$
we set $g$ to be the weak limit of
$\int_{X}F(x)\Phi(x)d\mu(x)$ in $\mathcal{H}$. It then follows that $F(x)=V_\Psi f(x)$, where $f:=C_{\Psi,\Phi}^{-1}g$, since
\begin{align*}
 V_\Psi f(x) &=\langle C_{\Psi,\Phi}^{-1}g,\Psi(x)\rangle=\langle g,(C_{\Psi,\Phi}^{-1})^\ast\Psi(x)\rangle\\
 &=\int_XF(y)\langle \Phi(y),(C_{\Psi,\Phi}^{-1})^\ast\Psi(x)\rangle d\mu(y)=R(F)(x)=F(x)
\end{align*}
It remains to show that if $F\in L^1(X,\mu)\cap L^\infty(X,\mu)$ it follows that the integral $\int_{X}F(x)\Phi(x)d\mu(x)$ converges weakly.
Let $h\in\mathcal{H}$
\begin{align*}
\big|\big\langle \int_{X}F(x)\Phi(x)d\mu(x),h\big\rangle\big| &\leq\int_{X}|F(x)||\langle\Phi(x),h\rangle |d\mu(x)\\
 &\leq \sup\limits_{x\in X}\|\Phi(x)\|_\mathcal{H}\ \|F\|_{L^1(X,\mu)}\ \|h\|_\mathcal{H}
\end{align*}
and hence by Riesz representation theorem $F\in dom(V_\Phi^\ast)$.
\hfill $\Box$

\vspace{0.3cm}
Observe that, unlike in the frame setting, there may exist $f\in\mathcal{H}$, s.t.
$V_\Psi f \notin L^2(X,\mu)$, which we will see later on in an example.
However, it can be easily seen that the images of $V_\Psi$ and $V_\Phi$ are subspaces of mutually dual spaces with the duality pairing 
$\langle F,H\rangle=\int_X F(x)\overline{H}(x)d\mu(x)$.

When dealing with reproducing pairs and continuous frames generated by a particular structure, for example by the action of a group 
representation  to a single window, three questions naturally arise:
Are there equivalent or sufficient conditions, independent of $f\in\mathcal{H}$, to obtain reproducing pairs or continuous frames? 
What can be said about the structure of the frame operator (resp. resolution operator)?
Given the mapping  $\Psi$ with a particular structure, is there another mapping $\Phi$ generated by the same structure such that
$(\Psi,\Phi)$ is a reproducing pair and the resolution operator is the identity operator?

In \cite{grmopa86} the authors gave a sufficient answer to these questions for the special case that $X$ is a locally compact group and
$\mu$ the left Haar measure. If $\pi:G\rightarrow\mathcal{H}$ is a square-integrable group representation, i.e.  if it is irreducible and
\begin{equation*}
\mathcal{A}:=\Big\{\psi\in\mathcal{H}:\ \int_G\left|\langle \pi(x)\psi,\psi\rangle\right|^2d\mu(x)<\infty\Big\}\neq\{0\}
\end{equation*}
then there exists a unique self-adjoint operator $L$ with domain $\mathcal{A}$, s.t. for all
$\psi,\varphi\in\mathcal{A}$ the following orthogonality relation holds
\begin{equation*}
 \int_G\langle f_1,\pi(x)\psi\rangle\overline{\langle f_2,\pi(x)\varphi\rangle}d\mu(x)=\langle L\varphi, L\psi\rangle\langle f_1,f_2\rangle
\end{equation*}
Elements of $\mathcal{A}$ are called admissible windows. If $G$ is unimodular, then $L$ is a multiple of the identity. Hence, we see that the
resolution operator is the identity after normalization. 

Regarding only systems arising from square-integrable group representations is nevertheless rather restrictive. This is why
we will introduce more flexible transforms in the next section.\\ \\
Within this paragraph we list some important results on LCA groups and its Fourier 
analysis. For a thorough introduction, see the standard text books \cite{fo95,lo53}. The most fundamental examples of 
LCA groups in harmonic analysis are the additive groups $\mathbb{R},\ \mathbb{Z},\ \mathbb{R}/\mathbb{Z}$ and $\mathbb{Z}/ N\mathbb{Z}$ and 
their d-fold products. Their relation is depicted in the following diagram.
\begin{center}
\begin{tikzpicture}[scale=0.7]
\draw[->,thick] (0,1.8)--(0,1)node[pos=0.5,left,font=\footnotesize]{Sampling};
\draw[->,thick] (0.5,0.5)--(3.3,0.5)node[pos=0.45,below,font=\footnotesize]{Periodization};
\draw[->,thick] (0.5,2.3)--(3.3,2.3)node[pos=0.45,above,font=\footnotesize]{Periodization};
\draw[->,thick] (3.9,1.8)--(3.9,1)node[pos=0.5,right,font=\footnotesize]{Sampling};
\draw (0,0.5) node {$\mathbb{Z}$};
\draw (0,2.3) node {$\mathbb{R}$};
\draw (3.9,2.3) node {$\mathbb{R}/\mathbb{Z}$};
\draw (4.1,0.5) node {$\mathbb{Z}/N\mathbb{Z}$};
\end{tikzpicture}
\end{center}
Every LCA group possesses a unique translation invariant 
measure on $G$ (up to a constant factor) called the Haar measure, denoted by $dx$. Convolution of two functions is given by
$f\ast g(y):=\int_G f(x)g(x^{-1}y)dx$.
It follows by Riesz-Thorin theorem that convolution with a fixed function $g\in L^1(G)$
is a bounded operator in $L^p(G)$, and $\|f\ast g\|_p\leq \|f\|_p\|g\|_1$, for $1\leq p\leq\infty$.

A character $\xi$ is a continuous homomorphism from $G$ to the torus $\mathbb{T}$, i.e. $\xi(xy)=\xi(x)\xi(y)$ and $|\xi(x)|=1$. 
The dual group $\widehat G$ of $G$ is the set of all characters of $G$. It is an LCA group with pointwise multiplication and the 
topology of compact convergence on $G$. 
The Pontryagin duality theorem states that any LCA group is ``reflexive'', i.e. the dual group of $\widehat G$ is isomorphic 
to $G$. The dual groups of the fundamental examples are given by
$\widehat{\mathbb{R}}\cong\mathbb{R},\  \widehat{\mathbb{R}/\mathbb{Z}}\cong\mathbb{Z},\ 
 \widehat{\mathbb{Z}}\cong\mathbb{R}/\mathbb{Z}$ and  $\widehat{\mathbb{Z}/N\mathbb{Z}}\cong\mathbb{Z}/N\mathbb{Z}$.
 
Now we are able to define the Fourier transform on $L^1(G)$ by
\begin{equation*}
 \hat f(\xi):=\int_G \overline{\xi(x)}f(x)dx,\ \ \xi\in\widehat G
\end{equation*}
It can be shown that this definition extends to an isometric isomorphism from $L^2(G)$ to $L^2(\widehat G)$ if the Haar
measure on $\widehat G$ is appropriately normalized, i.e. 
$\|f\|_2=\|\hat f\|_2$ and that Parseval's formula holds, i.e. $\langle f,g\rangle=\langle \hat f,\hat g\rangle,\ \forall f,g\in L^2(G)$. 
In addition, if $f,g\in L^2(G)$ and $f\ast g\in L^2(G)$, it follows that $(f\ast g)\ \widehat{}\ (\xi)=\hat f(\xi) \hat g(\xi)$.

\section{The continuous nonstationary Gabor transform on LCA groups}\label{cont-nonstationary gabor}
With the three questions from the previous section in mind we will now
consider continuous systems on $L^2(G)$ motivated by nonstationary Gabor frames. These discrete systems were introduced by 
Balazs et al. in \cite{nsdgt10} in order to gain more flexibility in analyzing signals with specific time-frequency characteristics. 
A major advantage of nonstationary Gabor frames is that while it allows adaptivity, it still guarantees perfect reconstruction, i.e.
resolution of the identity. Even more, under certain conditions, for a family $\{\psi_n\}_{n\in\mathbb{Z}}$ of compactly supported 
(resp. band-limited) windows the discrete frame operator is diagonal (resp. diagonal in the Fourier domain) and therefore allows for easy, 
and consequently fast, inversion.
 This property has been first studied by Daubechies et al. in 
\cite{daubgromay86} and is called the ``painless case'' in their study.

\subsection{Translation invariant systems}

Throughout the rest of this paper we assume $\mathcal{H}=L^2(G)$, where $G$ is a second countable LCA group. In particular this assumption
implies that $L^2(G)$ is separable and both $G$ and $\widehat G$ are $\sigma$-compact. As a consequence, it follows that the Haar measures 
on $G$ and $\widehat G$ are $\sigma$-finite.
The translation operator on $G$ is given by $T_z f(x):=f(z^{-1}x),\ x,z\in G$ and its Fourier transform 
is $\widehat{T_z f}(\xi)=\xi(z^{-1})\hat f(\xi),\ \xi\in\widehat G$.
Now let $\psi_y, \varphi_y\in L^2(G)$, for all $y\in Y$, where $(Y,\mu)$ is a measure space
with $\sigma$-finite measure $\mu$. For  $(x,y)\in G\times Y$ we define
\begin{equation*}
 \Psi(x,y):=T_{x}\psi_y\hspace{0.5cm} and \hspace{0.5cm}  \Phi(x,y):=T_{x}\varphi_y
 \end{equation*}
 and the continuous nonstationary Gabor transform (CNSGT) by
 \begin{equation*}
 V_\Psi f(x,y):=\langle f,\Psi(x,y)\rangle 
 \end{equation*}
 Then the following result on the structure of the reproducing (resp. frame) operator holds.

\begin{theorem}\label{cnsgt}
If there exist $A, B, C>0$, s.t.
 \begin{equation}\label{rep-pair-con}
  A\leq |m_{\Psi,\Phi}(\xi)|\leq B,\ for\ a.e.\ \xi\in\widehat{G}
 \end{equation} 
 where
  \begin{equation}
m_{\Psi,\Phi}(\xi):=\int_Y \overline{\widehat{\psi_y}(\xi)}\widehat{\varphi_y}(\xi)d\mu(y)
 \end{equation}
 and
   \begin{equation}\label{L1-con}
\int_Y \big|\widehat{\psi_y}(\xi)\widehat{\varphi_y}(\xi)\big|d\mu(y)\leq C,\ for\ a.e.\ \xi\in\widehat{G}
 \end{equation}
 then $(\Psi,\Phi)$ is a 
 reproducing pair for $L^2(G)$.
 The resolution operator is then given weakly by
 \begin{equation}\label{res-op-rep}
  C_{\Psi,\Phi} f=\mathcal{F}^{-1}(m_{\Psi,\Phi}\cdot\mathcal{F}(f))
 \end{equation}
 If $\Psi=\Phi$, then $\Psi$ is Bessel if and only if the upper bound in (\ref{rep-pair-con}) is satisfied and a continuous frame with 
 frame operator $S_{\Psi}=C_{\Psi,\Psi}$ and frame bounds $A,B$ if and only if condition (\ref{rep-pair-con}) is satisfied. In particular, 
 the frame is tight if condition (\ref{rep-pair-con}) becomes an equality.
\end{theorem}
\textbf{Proof:} Let $f_1,f_2\in 
L^1(G)\cap L^2(G)$, $\psi_y,\varphi_y\in L^2(G)$ and assume that (\ref{rep-pair-con}) and (\ref{L1-con}) hold.
Observe that $\langle f,T_{x}\psi_y\rangle=f\ast \psi_y^\ast(x)$, where $g^\ast(x):=\overline{g}(x^{-1})$ is the involution of $g$. 
Since $f\in L^1(G)$ it follows that $f\ast \psi_y^\ast\in L^2(G)$ and therefore $(f\ast \psi_y^\ast)\ \widehat{}\ (\xi)=\hat f(\xi)
\overline{\widehat \psi_y(\xi)}$.
Using this consideration and Parseval's formula we get
\begin{align*}
 \langle C_{\Psi,\Phi} f_1,f_2\rangle&=\int_Y\int_G\langle f_1,T_{x}\psi_y\rangle\overline{\langle f_2,T_{x}
 \varphi_y\rangle} dxd\mu(y)\\
&=\int_Y\int_{\widehat{G}}\hat{f}_1(\xi)\overline{\hat{f}}_2(\xi)
\overline{\widehat{\psi_y}(\xi)}\widehat{\varphi_y}(\xi)d\xi d\mu(y)\\
 &=\int_{\widehat G}m_{\Psi,\Phi}(\xi)\hat{f}_1(\xi)\overline{\hat{f}}_2(\xi)d\xi\\ &=\langle \mathcal{F}^{-1}(m_{\Psi,\Phi}\cdot\mathcal{F}
 (f_1)),f_2\rangle
\end{align*}
where condition (\ref{L1-con}) guarantees that Fubini's theorem is applicable. With the usual density argument, 
$C_{\Psi,\Phi}$  extends to a continuous operator on $L^2(G)$ which is bounded with bounded 
inverse by (\ref{rep-pair-con}).

It remains to show that if $\Psi$ is Bessel, it follows that the frame operator is given by (\ref{res-op-rep}). 
By previous calculation we get
\begin{equation*}
 \langle S_{\Psi} f,f\rangle=\int_Y\int_{\widehat G}|\hat f(\xi)|^2
|\widehat{\psi_y}(\xi)|^2 d\xi d\mu(y)\leq B\|f\|^2_2
\end{equation*}
Consequently, Fubini's theorem is again applicable and the frame operator is given by $S_\Psi f=\mathcal{F}^{-1}(m_{\Psi}\cdot
\mathcal{F}(f))$. It is easy to see that $S_\Psi$ is bounded with bounded inverse only if the symbol
$m_\Psi$ is essentially bounded from above and below.\hfill $\Box$

\vspace{0.3cm}
With a slight misuse of terminology we call $\{\psi_y\}_{y\in Y}$ admissible if (\ref{rep-pair-con}) is satisfied and 
 $\{(\psi_y,\varphi_y)\}_{y\in Y}$ cross-admissible, if conditions 
(\ref{rep-pair-con}) and (\ref{L1-con}) are satisfied.
Note that the inverse of a Fourier multiplier is given 
by another Fourier multiplier with the inverse symbol, i.e. $C_{\Psi,\Phi}^{-1}f=\mathcal{F}^{-1}(m_{\Psi,\Phi}^{-1}\cdot\mathcal{F}(f))$.

\begin{corollary}\label{cor-trans-inv}
If both $\Psi$ and $\Phi$ are Bessel, then $(\Psi,\Phi)$ is a reproducing system if and only if there exists $A>0$, s.t. 
$A\leq|m_{\Psi,\Phi}(\xi)|$, for a.e. $\xi\in\widehat G$ and  $C_{\Psi,\Phi}$ is given by (\ref{res-op-rep}).\\
Suppose that $C_{\Psi,\Phi}$ is given by (\ref{res-op-rep}). The resolution operator is the identity 
in $L^2(G)$ if and only if $m_{\Psi,\Phi}(\xi)=1$, for a.e. $\xi\in\widehat G$.\\
The canonical dual of a translation invariant frame $\Psi$ is another translation invariant system $\Phi(x,y):=T_xS_\Psi^{-1}\psi_y$
\end{corollary}
\textbf{Proof:} Let $\xi\in\widehat G$ s.t.  $m_{\Psi}(\xi)\leq B_\Psi$ and  $m_{\Phi}(\xi)\leq B_\Phi$, then it follows by Cauchy-Schwarz
inequality
\begin{equation*}
 |m_{\Psi,\Phi}(\xi)|\leq \int_Y|\hat \psi_y(\xi)\hat\varphi_y(\xi)|d\mu(y)\leq \left(m_\Psi(\xi) m_\Phi(\xi)\right)^{1/2}\leq 
 \left(B_\Psi B_\Phi\right)^{1/2}
\end{equation*}
Observe that this bound holds for a.e. $\xi\in\widehat G$, since the set in $\widehat G$ where the bound is violated is a union of two null 
sets. Consequently, (\ref{L1-con}) holds and Fubini's theorem is applicable. Finally, $A\leq|m_{\Psi,\Phi}(\xi)|$, for a.e. $\xi\in\widehat
G$, guarantees that $C_{\Psi,\Phi}$ is continuously invertible.

Since the Fourier transform is an onto isometry it follows that the equation  $f=\mathcal{F}^{-1}(m_{\Psi,\Phi}^{-1}\cdot\mathcal{F}(f))$
 holds for all $f\in L^2(G)$ if and only if the Fourier transform of $f$ is not altered in $L^2(G)$, i.e. $m_{\Psi,\Phi}(\xi)=1$, for 
 a.e. $\xi\in\widehat G$.
 
 To proof the last assertion we only have to show that the translation operator commutes with the inverse frame operator. But this is 
 obviously the case since the inverse frame operator is a Fourier multiplier operator and translation corresponds to character 
 multiplication in Fourier domain.
\hfill $\Box$

\begin{remark}
 This result shows that if $\Psi,\ \Phi$ are both Bessel, then $(\Psi,\Phi)$ is a reproducing pair if the functions
 $\psi_y(\xi),\varphi_y(\xi)$ are not orthogonal, or almost orthogonal, in $L^2(Y,\mu)$ for almost every $\xi\in\widehat G$.
\end{remark}

\subsection{Character invariant systems}

Now we multiply the windows $\psi_y$ with a character $\xi\in\widehat G$ instead of translating them, i.e. we consider the operator
$M_\xi f(x):=\xi(x)f(x)$ and the mappings 
\begin{equation*}
 \Psi(\xi,y):=M_{\xi}\psi_y\hspace{0.5cm}  and \hspace{0.5cm}   \Phi(\xi,y):=M_{\xi}\varphi_y
 \end{equation*}
 where  $(\xi,y)\in\widehat{G}\times Y$ and derive a similar result as in Theorem \ref{cnsgt}.
\begin{corollary} 
The pair of mappings $(\Psi,\Phi)$ is a reproducing pair for  $L^2(G)$, if there exist 
 $A, B, C>0$ s.t.
 \begin{equation}\label{sym-con2}
  A\leq |m_{\Psi,\Phi}(x)|\leq B,\ for\ a.e.\ x\in G
 \end{equation} 
 where
 \begin{equation}
 m_{\Psi,\Phi}(x):=\int_Y \overline{\psi_y(x)}\varphi_y(x)d\mu(y)
 \end{equation}
  and
   \begin{equation}
\int_Y \big|\psi_y(x)\varphi_y(x)\big|d\mu(y)\leq C,\ for\ a.e.\ x\in G
 \end{equation}
 The resolution operator is weakly given by
 \begin{equation}
  C_{\Psi,\Phi} f=m_{\Psi,\Phi}\cdot f
 \end{equation}
 If $\Psi=\Phi$, then $\Psi$ is Bessel if and only if the upper bound in (\ref{sym-con2}) is satisfied and a continuous frame with 
 frame operator $S_{\Psi}=C_{\Psi,\Psi}$ and frame bounds $A,B$ if and only if condition (\ref{sym-con2}) is satisfied. In particular, 
 the frame is tight if condition (\ref{sym-con2}) becomes an equality.
\end{corollary}
\textbf{Proof:} If one uses that $\widehat{M_\omega f}(\xi)=T_\omega\hat f(\xi)$ which implies 
$\langle f,M_{\xi} \psi_y\rangle=\langle \hat{f},T_{\xi} \widehat{\psi_y}\rangle$ 
and $\mathcal{F}_{\widehat{G}}\mathcal{F}_Gf(x)=f(x^{-1})$ one can follow the proof of Theorem \ref{cnsgt} step by step.\hfill $\Box$

\begin{example}\label{NSGT-examples}
Let us apply these results to two short examples with $G=(\mathbb{R},+)$. Observe that in this situation $\widehat G\cong G$.
For the short-time Fourier system
$\Psi(x,\omega)=M_{\omega}T_x\psi,\ \Phi(x,\omega)=M_{\omega}T_x\varphi,\ (x,\omega)\in \mathbb{R}^{2d}$ with $\mu$ the Lebesgue measure, 
one gets the Fourier symbol 
$m_{\Psi,\Phi}(\xi)=\langle \varphi,\psi\rangle$ and the well-known inversion 
formula
\begin{equation*}
 f=\frac{1}{\langle \varphi,\psi\rangle}\int_{\mathbb{R}^{2d}}\langle f,M_\omega T_x \psi\rangle M_\omega T_x\varphi dxd\omega
\end{equation*}
The second example is chosen to show that the theory also applies to discrete measure spaces with weighted counting measure. 
Consider the semi-discrete wavelet system with dyadic scale grid, i.e.
$\Psi(x,j)=T_xD_{2^j}\psi,\ (x,j)\in \mathbb{R}\times\mathbb{Z}$, with the dilation operator $D_af(x):=a^{-1/2}f(x/a)$ and $Y=\mathbb{Z}$ 
equipped with the weighted counting measure $\mu(j)=2^{-j}$. The Fourier symbol then reads
\begin{equation*}
 m_{\Psi}(\xi)=\sum_{j\in\mathbb{Z}}|\hat\psi(2^j\xi)|^2
\end{equation*}
It is not difficult to verify that, if $\hat\psi$ is continuous and the following two conditions hold,
\begin{enumerate}[(i)]
 \item $\exists\ \xi_0\neq 0,$ s.t. $\inf_{a\in[1,2]}|\hat\psi(a\xi_0)|>0$ 
 \item $\exists\ C>0$, s.t. $|\hat\psi(\xi)|^2\leq\frac{C|\xi|}{(1+|\xi|)^2},\ \forall \xi\in\mathbb{R}$
\end{enumerate}
 $m_{\Psi}$ is essentially bounded from below by $(i)$ and from above by $(ii)$. The canonical dual frame is another 
semi-discrete wavelet system. This can be seen if we use Corollary \ref{cor-trans-inv} and the observation that $m_{\Psi}(2^j\xi)=
m_{\Psi}(\xi)$ for all $j\in\mathbb{Z}$.
\begin{align*}
 S_\Psi^{-1} D_{2^j}\psi&=\mathcal{F}^{-1}\big(m_{\Psi}^{-1}D_{2^{-j}}\hat\psi\big)=\mathcal{F}^{-1}\big(m_{\Psi}^{-1}(2^j\ \cdot\ )
 D_{2^{-j}}\hat\psi\big)\\
 &=\mathcal{F}^{-1}\big(D_{2^{-j}}(m_{\Psi}^{-1}\hat\psi)\big) =D_{2^j}S_\Psi^{-1}\psi
 =:D_{2^j}\widetilde\psi
\end{align*}
\end{example}
 
 \begin{remark}
Most of the common continuous transforms used in signal processing can be written as a translation invariant system
 and can therefore be treated with the previous results. Besides the short-time Fourier transform and the continuous wavelet transform,
 those are for example
the continuous shearlet transform \cite{kula12} or the continuous curvelet transform \cite{cado05}, just to mention a few.\\
Observe that, although Theorem \ref{cnsgt} provides an equivalent condition for finite frame systems on $\mathbb{C}^N$ it does not supply 
a condition for discrete frames on $L^2(\mathbb{R}^d)$.
\end{remark}
 
 \section{Reproducing pairs for the affine Weyl-\\ Heisenberg group}\label{aWH-section}

Now we want to further reduce the level of abstractness.
Recently, greater effort has been put in the study of the affine Weyl-Heisenberg group $G_{aWH}$ and its representations, see 
\cite{daforastte08,hola95,kato93,torr1,torr2}, because it contains both the 
affine group and the Weyl-Heisenberg group as subgroups, the underlying groups of the continuous 
wavelet transform and the short-time Fourier transform. Consequently, there is a wide range of transforms arising from this group
 and its subsets.
 
Topologically
$G_{aWH}$ is isomorphic to $\mathbb{R}^{2d}\times\mathbb{R}^\ast\times\mathbb{T}$ with the group law given by
\begin{center}
 $(x,\omega,a,\tau)\cdot(x',\omega',a',\tau')=(x+ax',\omega+\omega'/a,aa',\tau\cdot\tau'\cdot e^{-2\pi i\omega'\cdot x/a})$
\end{center}
with the neutral element  $e=(0,0,1,1)$ and inverse element 
\begin{center}
$(-x/a,-a\omega,1/a,e^{2\pi i\omega\cdot x/a}/\tau)$
\end{center}
The affine Weyl-Heisenberg group is unimodular and the Haar measure is given by $d\mu(x,\omega,a,\tau)=dxd\omega|a|^{-1}dad\tau$.
A unitary representation of $G_{aWH}$ on $L^2(\mathbb{R}^d)$ is given by
\begin{center}
 $\pi(x,\omega,a,\tau)\psi=\tau M_\omega T_x D_a \psi$
\end{center}
where the basic time-frequency operators on $\mathbb{R}^d$ are given by
\begin{center}
$T_x f(t)=f(t-x),$\ \ $M_\omega f(t)=e^{2\pi i\omega t}f(t)$,\ \  $D_a f(t)=|a|^{-d/2}f(t/a)$
\end{center}
Since $G_{aWH}$ is a locally compact group, one is at first interested if this representation is square-integrable.
Unfortunately, this is not the case because, loosely speaking, the group is too big.

To overcome this obstacle Torr{\'e}sani \cite{torr1} suggested to regularize the Haar measure by multiplying it with a weight function 
$\rho(\omega)$ and showed that under certain conditions this also leads to tight continuous frames. A different approach in the same paper 
considered subgroups of the affine Weyl-Heisenberg group to obtain square-integrability. For example, if $d=1$, the section
$(x,\eta_\lambda(a),a,\tau)$ with $\eta_\lambda(a)=\lambda\left(\frac{1}{a}-1\right),\ \lambda\in\mathbb{R}$ is a
subgroup of $G_{aWH}$ and its representation is square-integrable
with left Haar measure $\frac{dxda}{|a|^2}$ and 
\begin{equation*}
 m_{\Psi,\Phi}\equiv\int_\mathbb{R}\overline{\hat\psi(\xi)}\hat\varphi(\xi)\frac{d\xi}{|\xi+\lambda|}
\end{equation*}
Within the scope of this paper we do not restrict to square-integrable representations but want to use the results from the previous 
section. The key to reproducing pairs or continuous frames lies in an appropriate restriction of the group parameters and the choice of
a measure $\mu$ on those subsets. Suppose that
$\beta:\mathbb{R}^n\rightarrow \mathbb{R},\ \eta:\mathbb{R}^n\rightarrow \mathbb{R}^d$, with $\beta,\eta$ piecewise continuous, 
$1\leq n\leq d$ and $\{\omega\in\mathbb{R}^n:\ \beta(\omega)=0\}$ is a null-set of $\mathbb{R}^n$. We
consider 
\begin{equation*}
G_{\beta,\eta}:=\Big\{(x,\eta(\omega),\beta(\omega),1): (x,\omega)\in \mathbb{R}^{d+n}\Big\}\subset G_{aWH}
\end{equation*}
together with the mapping
\begin{equation*}
\Psi(x,\omega)=M_{\eta(\omega)}T_x D_{\beta(\omega)}\psi,\ \ \psi\in L^2(\mathbb{R}^d)
\end{equation*}
The mapping $\Psi$ can be rewritten as $\Psi(x,\omega)=e^{2\pi i x\cdot\eta(\omega)}\widetilde\Psi(x,\omega)$, 
where $\widetilde\Psi(x,\omega)=T_x M_{\eta(\omega)}D_{\beta(\omega)}\psi$ is a translation invariant system.
By Lemma \ref{reprod-equiv-lemma} we can apply the recipe for translation invariant systems from Theorem \ref{cnsgt}.

As a measure on $G_{\beta,\eta}$ we take
$d\mu(x,\omega):=d\mu_s(x,\omega):=|\beta(\omega)|^{s-d}dxd\omega$, $s\in\mathbb{R}$. This particular definition of the measure 
$\mu_s$ is justified by two points. For one thing, the behavior of such a system is mainly depending on the scaling function $\beta$.
For another thing, it is necessary to introduce the parameter $s$ because altering $\beta,\eta$ may 
require a different choices of $s$ to ensure the existence of continuous frames. 
 To see this we assume that $d=1$ and consider on the one hand the 
continuous wavelet transform, i.e. $\eta\equiv 0$, $\beta(\omega)=\omega$. In this case, as 
$m_\Psi(\xi)=|\xi|^{-1-s}\int_\mathbb{R}|\hat\psi(a)|^2|a|^{s}da$,
there is no window $\psi\in L^2(\mathbb{R})$ such that this system forms a continuous frame if $s\neq-1$.
On the other hand for the setup $\eta(\omega)=\omega$, 
$\beta(\omega)=(1+|\omega|)^{-1}$, no 
continuous frame exists if $s\neq 1$, which we will explain more in more detail in Example \ref{ex-alpha=1}.
These two short examples also indicate that, in most cases, there is no freedom in the choice of the parameter $s$ to obtain a frame.

Theorem \ref{cnsgt} gives the following two sufficient conditions for $(\Psi,\Phi)$ to form a reproducing system for $L^2(\mathbb{R}^d)$
\begin{equation*}
 A\leq |m_{\Psi,\Phi}(\xi)|\leq B,\ for\ a.e.\ \xi\in\mathbb{R}^d
\end{equation*}
where
\begin{equation*}
m_{\Psi,\Phi}(\xi)=\int_{\mathbb{R}^n}\overline{\hat\psi\big(\beta(\omega)(\xi-\eta(\omega))\big)}\hat\varphi\big(\beta(\omega)
(\xi-\eta(\omega))\big)|\beta(\omega)|^{s}d\omega
\end{equation*}
and 
\begin{equation*}
\int_{\mathbb{R}^n}\Big|\hat\psi\big(\beta(\omega)(\xi-\eta(\omega))\big)\hat\varphi\big(\beta(\omega)(\xi-\eta(\omega))\big)\Big|
 |\beta(\omega)|^{s}d\omega<\infty
\end{equation*}
This result can be found in different articles where special attention is given to particular choices of $\beta,\eta$.
In \cite{hola95} composite frames are introduced whereas in \cite{daforastte08,dalo06} a special focus is on the
$\alpha$-transform and its uncertainty principles.
By choosing $\eta(\omega)=\omega$ we get a transform whose time-frequency resolution is frequency dependent. 
This is of particular interest for example in audio processing if one wants to construct a transform following the time-frequency 
resolution of the human auditory system, see \cite{neccxxl13}.

 \begin{example}\label{ex-alpha=1} The following example has been introduced in \cite{hola95} as composite frames.
 Let $d=n=1$, and $\eta(\omega)=\omega$, $\beta(\omega)=(1+|\omega|)^{-1}$. Substituting
 $z=\beta(\omega)(\xi-\omega)$ yields
 \begin{equation*}
 m_{\Psi,\Phi}(\xi)=\int_\mathbb{R}\overline{\hat\psi\Big(\frac{\xi-\omega}{1+|\omega|}\Big)}\hat\varphi\Big(\frac{\xi-\omega}{1+|\omega|}
 \Big)\frac{d\omega}{(1+|\omega|)^s}
  \end{equation*}
  \begin{equation*}
  =|1+\xi|^{1-s}\int_{-1}^\xi\overline{\hat\psi(z)}\hat\varphi(z)\frac{dz}{(1+z)^{2-s}}+
  |1-\xi|^{1-s}\int_{\xi}^1\overline{\hat\psi(z)}\hat\varphi(z)\frac{dz}{(1-z)^{2-s}}
 \end{equation*}
 It is easy to see that $ m_{\Psi,\Phi}$ fails to have either upper or lower bound if $s\neq1$ and $ m_{\Psi,\Phi}$ reads
  \begin{equation}\label{m-alpha-1}
 m_{\Psi,\Phi}(\xi)=\int_{-1}^\xi\overline{\hat\psi(z)}\hat\varphi(z)\frac{dz}{1+z}+\int_{\xi}^1\overline{\hat\psi(z)}\hat\varphi(z)
 \frac{dz}{1-z}
\end{equation}
This expression allows for explicit calculation of $m_{\Psi,\Phi}$ for many choices of $\psi,\varphi$, take
 for example $\hat\psi(\xi)=(1-\xi)\chi_{A}(\xi)$ and $\hat\varphi(\xi)=(1+\xi)\chi_{A}(\xi)$, where $A:=[-1,1]$, then one gets
 $m_{\Psi,\Phi}(\xi)=2$, for $|\xi|>1$ and $m_{\Psi,\Phi}(\xi)=3-\xi^2$, for $|\xi|\leq 1$.
 \end{example}
 
 In the following we answer the question whether there always exists a ``dual system'' with the same structure and examine drawbacks
 and advantages of  reproducing pairs all with the aid of Example \ref{ex-alpha=1}.

 \begin{proposition}
For the transform of Example \ref{ex-alpha=1} there is no reproducing pair $(\Psi,\Phi)$, s.t. $\Psi$ and $\Phi$ are Bessel and 
$C_{\Psi,\Phi}=Id_{L^2(\mathbb{R})}$.
 \end{proposition}
\textbf{Proof:} Since $\Psi,\Phi$ are Bessel we get by Corollary \ref{cor-trans-inv} that $C_{\Psi,\Phi}=Id_{L^2(\mathbb{R})}$
if and only if $m_{\Psi,\Phi}(\xi)=1$, for a.e. $\xi\in\mathbb{R}$. W.l.o.g. we assume that $\hat\psi,\hat\varphi$ are real-valued 
functions. Otherwise use $Re(m_{\Psi,\Phi})$ instead of $m_{\Psi,\Phi}$ for the following arguments. In order to obtain a contradiction 
we assume that there exist  $\psi,\varphi\in L^2(\mathbb{R})$, s.t. $m_{\Psi,\Phi}=1$, for a.e. $\xi\in\mathbb{R}$. Then 
$m_{\Psi,\Phi}(\xi)=1$, for all $\xi\in[-1,1]$, since both summands in (\ref{m-alpha-1}) are continuous. Hence, it follows that 
$m_{\Psi,\Phi}'(\xi)=0$, for every  $\xi\in(-1,1)$. On the other hand, Lebesgue's differentiation theorem states that, for a.e. 
$\xi\in(-1,1)$, the derivative of  $m_{\Psi,\Phi}$ is given by
  \begin{equation*}
  m_{\Psi,\Phi}'(\xi)=\hat\psi(\xi)\hat\varphi(\xi)\left[\frac{1}{1+\xi}-\frac{1}{1-\xi}\right]=-\frac{2\xi}{1-\xi^2}
  \hat\psi(\xi)\hat\varphi(\xi)
  \end{equation*}
which implies
\begin{equation*}
\hat\psi(\xi)\hat\varphi(\xi)=0,\ for\ a.e.\ \xi\in(-1,1)
\end{equation*}
This finally yields the contradiction $m_{\Psi,\Phi}(\xi)=0,\ \forall\xi\in[-1,1]$.\hfill $\Box$

\begin{proposition}
In general it holds that if $(\Psi,\Phi)$ is a reproducing pair, neither $\Psi$ nor $\Phi$ has to be Bessel.
\end{proposition}
\textbf{Proof:} Consider the pair 
$\hat\psi(\xi)=(1-\xi)\chi_{A}(\xi)$, $\hat\varphi(\xi)=(1+\xi)\chi_{A}(\xi)$ from Example \ref{ex-alpha=1}. The
Fourier symbol  $m_{\Psi,\Phi}$ is bounded from above and below but neither $\Psi$ nor $\Phi$ are Bessel 
since $m_{\Psi}$ and $m_{\Phi}$ are unbounded. \hfill $\Box$

\vspace{0.3cm}
Finally, we show that the orbit of $V_\Psi$ possibly contains elements $V_\Psi f\notin L^2(X,\mu)$.
 Again, we consider  Example \ref{ex-alpha=1} with $\hat\psi(\xi)=(1+\xi)\chi_{A}(\xi)$ and $f \in L^1(\mathbb{R})\cap L^2(\mathbb{R})$, 
$|\hat f(\xi)|\geq 1,\ \forall \xi\in A$. For every $(\xi,\omega)\in A\times\mathbb{R}$ it holds $\beta(\omega)(\xi-\omega)\in A$.
Hence, it follows
\begin{align*}
 \|V_\psi f\|_{L^2(\mathbb{R}^2,\mu)}  &= \int_\mathbb{R}\int_\mathbb{R}|\hat f(\xi)|^2|\hat\psi(\beta(\omega)(\xi-\omega))|^2
\beta(\omega)d\xi d\omega\\
&\geq\int_\mathbb{R}\int_A|\hat\psi(\beta(\omega)(\xi-\omega))|^2\beta(\omega)d\xi d\omega\\
&=\int_\mathbb{R}\int_A\Big(1+\beta(\omega)(\xi-\omega))\Big)^2\beta(\omega)d\xi d\omega\\
&=\int_\mathbb{R}\int_A\big(1+2|\omega|\chi_{(-\infty,0]}(\omega)+\xi\big)^2\beta(\omega)^3d\xi d\omega\\
&=\int_\mathbb{R}\Big[C_1+\Big(C_2+2|\omega|\chi_{(-\infty,0]}(\omega)\Big)^2\ \Big]
\beta(\omega)^3 d\omega\\ &=\infty
\end{align*}

\section{Outlook and discussions}
It seems interesting to study discretization schemes for the CNSGT when starting from a semi-discrete system, see Example 
\ref{NSGT-examples}. Clearly, coorbit theory \cite{fora05}
could be applied here but this approach neither exploits that the index set $Y$ is already discrete, nor
the abelian group structure of $G$.

Using the results about frame multipliers, the more general concept of reproducing pairs might also be interesting for discrete systems.

Furthermore, as mentioned after Proposition \ref{sec:repordkern1}, a full characterization of the orbit of $V_\Psi$ is desirable. 
To this end, it might be worthwhile to construct and investigate Gelfand triples of those spaces, similar to the approach in 
\cite{jpaxxl09}. Moreover, as the images of $V_\Psi$ and $V_\Phi$ are mutually dual, an investigation in the context of partial 
inner product spaces \cite{antoine2009partial} seems appropriate.

\section*{Acknowledgement}
The authors like to acknowledge the financial support by the Austrian Science Fund (FWF) through the START-project 
FLAME (Frames and Linear Operators for Acoustical Modeling and Parameter Estimation): Y 551-N13.

\bibliographystyle{plain}
\bibliography{paperbib}

\begin{thebibliography}{10}

\bibitem{alanga93}
S.~T. Ali, J.~P. Antoine, and J.~P. Gazeau.
\newblock Continuous frames in {H}ilbert spaces.
\newblock {\em Ann. Phys.}, 222:1--37, 1993.

\bibitem{alanga93a}
S.~T. Ali, J.~P. Antoine, and J.~P. Gazeau.
\newblock Relativistic quantum frames.
\newblock {\em Ann. Phys.}, 222:38--88, 1993.

\bibitem{alanga00}
S.~T. Ali, J.~P. Antoine, and J.~P. Gazeau.
\newblock {\em Coherent states, wavelets and their generalizations}.
\newblock Graduate texts in contemporary physics. Springer, 2000.

\bibitem{jpaxxl09}
Jean-Pierre Antoine and Peter Balazs.
\newblock Frames and semi-frames.
\newblock {\em J. Phys. A: Math. Theor.}, 44:205201, 2011.

\bibitem{antoine2009partial}
J.P. Antoine and C.~Trapani.
\newblock {\em Partial Inner Product Spaces: Theory and Applications}.
\newblock Lecture notes in mathematics. Springer-Verlag Berlin Heidelberg,
  2009.

\bibitem{nsdgt10}
P.~Balazs, M.~D\"orfler, N.~Holighaus, F.~Jaillet, and G.~Velasco.
\newblock Theory, implementation and applications of nonstationary {G}abor
  frames.
\newblock {\em J. Comp. Appl. Math.}, 236(6):1481--1496, 2011.

\bibitem{badokowto13}
P.~Balazs, M.~D{\"o}rfler, M.~Kowalski, and B.~Torr{\'e}sani.
\newblock Adapted and adaptive linear time-frequency representations: a
  synthesis point of view.
\newblock {\em IEEE Signal Processing Magazine (special issue: Time-Frequency
  Analysis and Applications)}, 30(6):20--31, 2013.

\bibitem{xxlmult1}
Peter Balazs.
\newblock Basic definition and properties of {B}essel multipliers.
\newblock {\em J. Math. Anal. Appl.}, 325(1):571--585, January 2007.

\bibitem{xxlbayasg11}
Peter Balazs, Dominik Bayer, and Asghar Rahimi.
\newblock Multipliers for continuous frames in {H}ilbert spaces.
\newblock {\em J. Phys. A: Math. Theor.}, 45:244023, 2012.

\bibitem{cado05}
E.~J. Cand{\`e}s and D.~L. Donoho.
\newblock Continuous curvelet transform: I. {R}esolution of the wavefront set.
\newblock {\em Appl. Comp. Harmon. Anal.}, 19:162--197, 2005.

\bibitem{raey}
Peter~G. Casazza, Gitta Kutyniok, and Mark~C. Lammers.
\newblock Duality principles in frame theory.
\newblock {\em Journal of Fourier Analysis and Applications}, 10(4):383--408,
  2004.

\bibitem{christ1}
O.~Christensen.
\newblock {\em An {I}ntroduction to {F}rames and {R}iesz {B}ases.}
\newblock {A}pplied and {N}umerical {H}armonic {A}nalysis. {B}irkh{\"a}user,
  2003.

\bibitem{Christensen2014198}
Ole Christensen and Say~Song Goh.
\newblock From dual pairs of gabor frames to dual pairs of wavelet frames and
  vice versa.
\newblock {\em Applied and Computational Harmonic Analysis}, 36(2):198 -- 214,
  2014.

\bibitem{daforastte08}
S.~Dahlke, M.~Fornasier, H.~Rauhut, G.~Steidl, and G.~Teschke.
\newblock Generalized coorbit theory, {B}anach frames, and the relation to
  $\alpha$-modulation spaces.
\newblock {\em {P}roc. {L}ondon {M}ath. {S}oc.}, 96(2):464--506, 2008.

\bibitem{dalo06}
S.~Dahlke, D.~Lorenz, P.~Maass, C.~Sagiv, and G.~Teschke.
\newblock The canonical coherent states associated with quotients of the affine
  {W}eyl-{H}eisenberg group.
\newblock {\em J. Appl. Funct. Anal.}, 3:215--232, 2008.

\bibitem{daubgromay86}
I.~Daubechies, A.~Grossmann, and Y.~Meyer.
\newblock Painless non-orthogonal expansions.
\newblock {\em J. Math. Phys.}, 27:1271--1283, 1986.

\bibitem{fo95}
G.~B. Folland.
\newblock {\em A course in abstract harmonic analysis}.
\newblock Studies in advanced mathematics. CRC Press, 1995.

\bibitem{fora05}
M.~Fornasier and H.~Rauhut.
\newblock Continuous frames, function spaces, and the discretization problem.
\newblock {\em {J}. {F}ourier {A}nal. {A}ppl. 11}, pages 244--287, 2005.

\bibitem{groe1}
K.~Gr{\"o}chenig.
\newblock {\em {F}oundations of {T}ime-{F}requency {A}nalysis}.
\newblock {A}ppl. {N}umer. {H}armon. {A}nal. {B}irkh{\"a}user {B}oston, 2001.

\bibitem{grmopa86}
A.~Grossmann, J.~Morlet, and T.~Paul.
\newblock Transforms associated to square integrable group representations.
\newblock {\em Ann. Inst. {H}enri {P}oincar{\'e}}, 45, 1986.

\bibitem{hola95}
J.~A. Hogan and J.~D. Lakey.
\newblock Extensions of the heisenberg group by dilations and frames.
\newblock {\em Appl. Comp. Harmon. Anal.}, 2:174--199, 1995.

\bibitem{jale14}
M.~S. Jakobsen and J.~Lemvig.
\newblock Reproducing formulas for generalized translation invariant systems on
  locally compact abelian groups.
\newblock submitted, 5 2014.

\bibitem{ka90}
G.~Kaiser.
\newblock {\em Quantum physics, relativity, and complex spacetime: towards a
  new synthesis}.
\newblock North Holland, 1990.

\bibitem{ka94}
G.~Kaiser.
\newblock {\em Generalized frames: Key to analysis and synthesis}, chapter~4,
  pages 78--98.
\newblock Birkh{\"a}user, 1994.

\bibitem{kato93}
C.~Kalisa and B.~Torr{\'e}sani.
\newblock $n$-dimensional affine {W}eyl-{H}eisenberg wavelets.
\newblock {\em {A}nn. {I}nst. {H}enri {P}oincar{\'e} {P}hys. {T}h'eor.},
  59(2):201--236, 1993.

\bibitem{kula12}
G.~Kutyniok and D.~Labate.
\newblock Introduction to shearlets.
\newblock In Kutyniok G. and D.~Labate, editors, {\em Shearlets: Multiscale
  Analysis for Multivariate Data}, Applied and Numerical Harmonic Analysis,
  chapter~1, pages 1--36. Birkh{\"a}user, 2012.

\bibitem{lo53}
L.~H. Loomis.
\newblock {\em An introduction to abstract harmonic analysis}.
\newblock D. van Nostrand Company, Inc., New York, 1953.

\bibitem{ma09}
S.~Mallat.
\newblock {\em A Wavelet Tour of Signal Processing - The Sparse Way}.
\newblock Academic Press, 3rd edition, 2009.

\bibitem{neccxxl13}
T.~Necciari, P.~Balazs, N.~Holighaus, and P.~S{\o}ndergaard.
\newblock The {ERB}let transform: An auditory-based time-frequency
  representation with perfect reconstruction.
\newblock In {\em ICASSP 2013}, 2013.

\bibitem{persampta13}
N.~Perraudin, N.~Holighaus, P.~Soendergaard, and P.~Balazs.
\newblock Gabor dual windows using convex optimization.
\newblock In {\em Proceeedings of the 10th International Conference on Sampling
  theory and Applications (SAMPTA 2013)}, 2013.

\bibitem{ranade06}
A.~Rahimi, A.~Najati, and Y.~N. Dehghan.
\newblock Continuous frames in {H}ilbert spaces.
\newblock {\em Methods Funct. Anal. Topol.}, 12(2):170--182, 2006.

\bibitem{stobalrep11}
Diana Stoeva and Peter Balazs.
\newblock Representation of the inverse of a frame multiplier.
\newblock submitted, -.

\bibitem{balsto09new}
Diana~T. Stoeva and Peter Balazs.
\newblock Invertibility of multipliers.
\newblock {\em Appl. Comp. Harmon. Anal.}, 33(2):292--299, 2012.

\bibitem{uncconv2011}
Diana~T. Stoeva and Peter Balazs.
\newblock Canonical forms of unconditionally convergent multipliers.
\newblock {\em J. Math. Anal. Appl.}, 399:252--259, 2013.

\bibitem{te01}
G.~Teschke.
\newblock {\em Waveletkonstruktion {\"u}ber {U}nsch{\"a}rferelationen und
  {A}nwendungen in der {S}ignalanalyse}.
\newblock PhD thesis, Universit{\"a}t Bremen, 2001.

\bibitem{torr1}
B.~Torr{\'e}sani.
\newblock Wavelets associated with representations of the {W}eyl-{H}eisenberg
  group.
\newblock {\em J. Math. Phys.}, 32(5), May 1991.

\bibitem{torr2}
B.~Torr{\'e}sani.
\newblock Time-frequency representations: wavelet packets and optimal
  decomposition.
\newblock {\em Annales de l'IU.H.P.}, (2):215--234, 1992.

\end{thebibliography}
\end{document}